\documentclass{article}
\usepackage{latexsym,amsmath,amssymb,amscd}
\begin{document}

\title{On a certain type of nonlinear hyperbolic equations derived from
astrophysical problems}
\author{Tetu Makino \footnote{Professor Emeritus at Yamaguchi University,
Japan / e-mail: makino@yamaguchi-u.ac.jp}}
\date{\today}
\maketitle

\newtheorem{Lemma}{Lemma}
\newtheorem{Proposition}{Proposition}
\newtheorem{Theorem}{Theorem}
\newtheorem{Definition}{Definition}
\newtheorem{Corollary}{Corollary}

\begin{abstract}
Investigations of spherically symmetric motions 
of self-gravitating gaseous stars governed by the
non-relativistic Newtonian gravitation theory or by the general
relativistic theory lead us to a certain type of non-linear
hyperbolic equations defined on a finite interval of the space variable.
The linearized principal part has regular singularities at the both ends of
the interval of the space variable. But the regularity loss caused by the 
singularities at the boundaries requires application of the Nash-Moser technique.
An abstract unified treatment of the problem is presented.\\



\end{abstract}

\section{Introduction}

We have investigated spherically symmetric motions of
gaseous stars governed by the non-relativistic Euler-Poisson equations
in \cite{OJM} and by the relativistic Einstein-Euler equations in
\cite{KJM}. See also \cite{EEdS} on the problem governed by
the Einstein-Euler-de Sitter equations with the
cosmological constant. In this article a unified formulation of the discussions
developed in the previous works is presented. This treatment using the Nash-Moser theorem is based on the work \cite{FE}.
Let us begin by giving an abstract settlement of the problem.\\

The set of equations we consider is
\begin{subequations}
\begin{eqnarray}
&&\frac{\partial y}{\partial t}-J(x,y,z)v=0, \label{Ea} \\
&&\frac{\partial v}{\partial t}+H_1(x,y,z,v)\mathcal{L}y+
H_2(x,y,z,v,w)=0 \label{Eb}
\end{eqnarray}
\end{subequations}
with

\begin{equation}
\mathcal{L}=-x(1-x)\frac{d^2}{dx^2}-
\Big(\frac{5}{2}(1-x)-\frac{N}{2}x\Big)\frac{d}{dx}+
\ell_1(x)x(1-x)\frac{d}{dx}+L_0(x). \label{2}
\end{equation}
Here $x$ runs on the interval $]0,1[$ and 
$\displaystyle z=x\frac{\partial y}{\partial x}, w=x\frac{\partial v}{\partial x}$.
$N$ is a parameter and we put the assumption \\

{\bf (B0):}  $N>4$.\\

Let us
denote by $\mathfrak{A}$ the set of all functions defined and
analytic on a neighborhood of $[0,1]$. 
Let us denote by $\mathfrak{A}^Q(U^p)$ the set of all analytic
functions $f(x, y_1, \cdots, y_p)$ of $x$ in a neighborhood 
of $[0,1]$ and $y_1 \cdots, y_p$ in a neighborhood $U$ of $0$ such that
$$f(x,y_1,\cdots,y_p)=\sum_{k_1+\cdots+k_p\geq Q }
a_{k_1\cdots k_p}(x)y_1^{k_1}\cdots y_p^{k_p}. $$

We put \\

{\bf (B1):} {\it We suppose $\ell_1, L_0 \in \mathfrak{A}$ and 
there is a neighborhood $U$ of $0$ such that $J\in \mathfrak{A}^0(U^2)$,  
$H_1\in \mathfrak{A}^0(U^3)$, 
$H_2\in \mathfrak{A}^2(U^4)$.}\\

Let us fix $T>0$ arbitrarily and fix functions
$y^*, v^* \in C^{\infty}([0,T]\times [0,1])$ such that
$y^*(t,x), z^*(t,x)=x\partial y^*/\partial x, v^*(t,x), 
w^*(t,x)=x\partial v^*/\partial x \in U$ for $0\leq \forall t \leq T,
0\leq\forall x \leq 1$.
We seek a solution $y, v \in C^{\infty}([0,T]\times [0,1])$ of
(\ref{Ea})(\ref{Eb}) of the form
\begin{equation}
y=y^*+\tilde{y},\qquad v=v^*+\tilde{v} \label{Var} 
\end{equation}
which satisfies
\begin{equation}
\tilde{y}|_{t=0}=0,\qquad \tilde{v}|_{t=0}=0. \label{Init}
\end{equation}

We suppose the following assumptions:\\

{\bf (B2):}  {\it We have
\begin{equation}
H_1(x,0,0,0)=J(x,0,0)^{-1}
\end{equation}
and there is a constant $C>1$ such that
\begin{equation}
\frac{1}{C}<J(x,0,0)<C
\end{equation}
for $\forall x \in U_0$.} \\

{\bf (B3):}  {\it We have
\begin{equation}
\partial_zJ\equiv 0,\quad
(\partial_zH_1)\mathcal{L}y+\partial_zH_2\equiv 0,
\quad \partial_wH_2\equiv 0 \label{B2}
\end{equation}
as $x\rightarrow 1$.}\\

Here a function $f$ defined and analytic on a neighborhood of
$(1,0,\cdots, 0)$ is said to satisfy $f\equiv 0$ as $x\rightarrow 1$ if
there is a function $\Omega$ defined and analytic on a neighborhood
of $(1,0,\cdots, 0)$ such that
$$f(x,y_1,\cdots, y_p)=(1-x)\Omega(x,y_1,\cdots, y_p).$$
Of course it is the case if 
$f(1, y_1,\cdots, y_p)=0$ for $\forall y_1, \cdots, \forall y_p$.
In the assumption {\bf (B3)} the functions in (\ref{B2}) are 
regarded as functions of $x, y, Dy, D^2y, v, Dv$. Here $D$ stands for $\partial/\partial x$. \\

Under the above situation, we are going to prove the 
following

\begin{Theorem}
There is a small positive number $\delta(T)$ and a large number $\mathfrak{K}$ such that, if
\begin{equation}
\max_{j+k\leq\mathfrak{K}}
\|\partial_t^j\partial_x^k(y^*, v^*)\|_{L^{\infty}}\leq\delta(T),
\end{equation}
there exists a solution $(y,v)$ of  (\ref{Ea})(\ref{Eb})(\ref{Var})(\ref{Init}).
\end{Theorem}

\section{Frame-work to apply the Nash-Moser(-Hamilton) theorem}

The equations for $\mathsf{w}:=(\tilde{y}, \tilde{v})^T$ turn out to be
\begin{subequations}
\begin{eqnarray}
&&\frac{\partial\tilde{y}}{\partial t}-J\tilde{v}-
(\Delta J)v^*=c_1, \label{Eqw_a} \\
&&\frac{\partial\tilde{v}}{\partial t}+
H_1\mathcal{L}\tilde{y}+
(\Delta H_1)\mathcal{L}y^*+\Delta H_2=c_2, \label{Eqw_b}
\end{eqnarray}
\end{subequations}
where
\begin{align*}
J&=J(x,y^*+\tilde{y},
z^*+\tilde{z}) \quad \mbox{with}\quad\tilde{z}=
x\frac{\partial\tilde{y}}{\partial x} \\
\Delta J&=J(x,y^*+\tilde{y}, z^*+\tilde{z})-
J(x,y^*, z^*), \\
c_1&=-\frac{\partial y^*}{\partial t}+J(x,y^*,z^*)v^*, \\
H_1&=H_1(x, y^*+\tilde{y},z^*+\tilde{z}, v^*+\tilde{v}, w^*+\tilde{w}) 
\quad\mbox{with}\quad \tilde{w}=x\frac{\partial\tilde{v}}{\partial x}, \\
\Delta H_1&=H_1(x,y^*+\tilde{y}, z^*+\tilde{z}, v^*+\tilde{v})-
H_1(x,y^*,z^*,v^*), \\
\Delta H_2&=H_2(x,y^*+\tilde{y}, z^*+\tilde{z}, v^*+\tilde{v},
w^*+\tilde{w}) 
-H_2(x,y^*,z^*, v^*, w^*), \\
c_2&=-\frac{\partial {v}^*}{\partial t}-
H_1(x,y^*, z^*, v^*)\mathcal{L}y^*-
H_2(x,y^*,z^*,v^*,w^*).
\end{align*}

We write the equations (\ref{Eqw_a})(\ref{Eqw_b}) as
\begin{equation}\mathfrak{P}(\mathsf{w})=\mathsf{c},
\end{equation}
 where 
$\mathsf{c}=(c_1,c_2)^T$. The domain of the nonlinear mapping
$\mathfrak{P}$ is $\mathfrak{U}$, the set of all functions
$\mathsf{w}=(\tilde{y},\tilde{v})^T \in \mathfrak{E}_0\times
\mathfrak{E}_0$ such that
\begin{equation}
|\tilde{y}|+|xD\tilde{y}|+|\tilde{v}|+
|xD\tilde{v}|<\epsilon_0, \label{U}
\end{equation}
where $D=\partial/\partial x$.
Here $\epsilon_0$ is so small that (\ref{U}) implies
$y(t,x)=y^*(t,x)+\tilde{y}(t,x), z=xDy, v, w=xDv \in U$ for $\forall t \in [0,T],
\forall x\in U_0$. We have defined
$\mathfrak{E}=C^{\infty}([0,T]\times[0,1])$ and
$\mathfrak{E}_0=\{\phi \in\mathfrak{E} | \  \phi|_{t=0}=0\}$.\\

We are going to apply the Nash-Moser(-Hamilton) theorem
(\cite[p. 171, III.1.1.1]{Hamilton}). \\

First we introduce gradings of norms on $\mathfrak{E}$
to make it a tame space in the sense of Hamilton \cite{Hamilton}.

 To do so,
we use a cut off function $\omega\in C^{\infty}(\mathbb{R})$ such that
$\omega(x)=1$ for $x\leq 1/3$, 
$0<\omega(x)<1$ for $1/3<x<2/3$ and
$\omega(x)=0$ for $2/3\leq x$. For a function $u$ of $x\in [0,1]$, we define
\begin{equation}
u^{[0]}(x):=\omega(x)u(x),\qquad
u^{[1]}(x):=(1-\omega(x))u(x).
\end{equation}

We consider the functional spaces
\begin{subequations}
\begin{eqnarray}
\mathfrak{E}_{[0]}&=&\{ u\in \mathfrak{E} |\  u=0 \quad\mbox{for}\quad 5/6\leq x\leq 1 \}, \\
\mathfrak{E}_{[1]}&=&\{ u\in\mathfrak{E} |\  u=0
\quad\mbox{for}\quad 0\leq x\leq 1/6\}
\end{eqnarray}
\end{subequations}
endowed with the gradings
$(\|\cdot|_{[\mu]n}^{(\infty)})_n,
(\|\cdot\|_{[\mu]n}^{(2)})_n$, $\mu=0,1$, defined as follows.

Put
$$\triangle_{[0]}:=x\frac{d^2}{dx^2}+\frac{5}{2}\frac{d}{dx}, $$
and
$$\triangle_{[1]}:=X\frac{d^2}{dX^2}+\frac{N}{2}\frac{d}{dX} \quad\mbox{with}\quad X=1-x,$$
and put
\begin{align*}
\|u\|_{[\mu]n}^{(\infty)}&=
\sup_{j+k\leq n}
\|(-\partial_t^2)^j(-\triangle_{[\mu]})^ku\|_{L^{\infty}}, \\
\|u\|_{[\mu]n}^{(2)}&=\Big(
\sum_{j+k\leq n}
\int_0^T\|(-\partial_t^2)^j
(-\triangle_{[\mu]})^ku\|_{[\mu]}^2dt\Big)^{1/2},
\end{align*}
where
\begin{align*}
\|u\|_{[0]}&=\Big(\int_0^1|u|^2x^{3/2}dx\Big)^{1/2}, \\
\|u\|_{[1]}&=\Big(\int_0^1|u|^2(1-x)^{N/2-1}dx\Big)^{1/2}.
\end{align*}

\begin{Proposition}
Let $\mu=0$ or $=1$. The gradings $(\|\cdot|_{[\mu]n}^{(\infty)})_n,
(\|\cdot\|_{[\mu]n}^{(2)})_n$ are equivalent and make $\mathfrak{E}_{[\mu]}$ a tame space.
\end{Proposition}

Proof. Let us consider $\mathfrak{E}_{[1]}$ and denote $X=1-x$ by $x$
and $\triangle_{[1]}$ by
$$\triangle=x\frac{d^2}{dx^2}+\frac{N}{2}\frac{d}{dx}.$$
( $\mathfrak{E}_{[0]}$ can be treated similarly.) Then $u \in \mathfrak{E}_{[1]}$ means $u \in C^{\infty}([0,1])$ and $u(x)=0$ for $x\geq 5/6$. 

Now, we can define the Fourier transformation
$Fu(\xi)$ of a function $u(x)$ of $x\in [0,+\infty)$ by
$$Fu(\xi):=\int_0^{\infty}K(\xi x)u(x)x^{N/2-1}dx, $$
where $K(X)$ is an entire function of $X\in\mathbb{C}$
given by
$$K(X)=2(\sqrt{X})^{-N/2+1}J_{N/2-1}(4\sqrt{X})=
2^{N/2}\sum_{k=0}^{\infty}
\frac{(-4X)^k}{k!\Gamma(N/2+k)},
$$
$J_{N/2-1}$ being the Bessel function of order $N/2-1$. Note that
$$F(-\triangle u)(\xi)=4\xi\cdot Fu(\xi)$$
and the inverse of $F$ is $F$ itself. Then we see that
$\mathfrak{E}_{[1]}$ endowed with $(\|\cdot\|_{[1]n}^{\infty})_n $
is a tame direct summand of the tame space
$$\mathfrak{F}=L_1^{\infty}
(\mathbb{R}\times[0,+\infty), d\tau\otimes\xi^{N/2-1}d\xi, \log(1+\tau^2+4\xi))$$ through the transformation
$$\mathcal{F}u(\tau, \xi)=
\frac{1}{\sqrt{2\pi}}\int e^{-\sqrt{-1}\tau t}Fu(t, \cdot)(\xi)dt $$
and its inverse applied to the space
$$\tilde{\mathfrak{E}}:=C_0^{\infty}((-2T, 2T)\times[0,2)),$$ into which 
the extension of functions of $\mathfrak{E}_{[1]}$
multiplied by a fixed cut off function $\chi(t)$ which vanishes for $t\leq -T$ can be imbedded,
and the space
$$\dot{\mathfrak{E}}:=\dot{C}^{\infty}(\mathbb{R}\times[0,+\infty))
:=\{ u |\  \forall j\forall k
\lim_{L\rightarrow \infty}\sup_{L\leq|t|, L\leq x}
|(-\partial_t^2)^j(-\triangle)^ku|=0\},$$
for which functions of $\mathfrak{E}_{[1]}$ are restrictions. In fact,
if we denote by $\mathfrak{e}:\mathfrak{E}_{[1]}\rightarrow
\tilde{\mathfrak{E}}$ the extension operator, and
by $\mathfrak{r}:\dot{\mathfrak{E}}\rightarrow \mathfrak{E}_{[1]}$
the restriction operator, then
the mappings $\mathcal{F}\mathfrak{e}:\mathfrak{E}_{[1]}\rightarrow
\mathfrak{F}$ and $\mathfrak{r}\mathcal{F}:\mathfrak{F}\rightarrow 
\mathfrak{E}_{[1]}$ are tame and the composition
$(\mathfrak{r}\mathcal{F})\circ(\mathcal{F}\mathfrak{e})$ is $\mbox{Id}_{\mathfrak{E}_{[1]}}$. For the details, see the proof of
\cite[p.137, II.1.3.6.]{Hamilton}. This implies that $\mathfrak{E}_{[1]}$ is tame 
with respect to $(\|\cdot\|_{[1]n}^{\infty)})_n$. On the other hand
$$\sqrt{\frac{N}{2}}\|u\|_{[1]}\leq \|u\|_{L^{\infty}}
\leq C\sup_{j\leq \sigma}\|(-\triangle)^ju\|_{[1]}$$
by the Sobolev imbedding theorem 
(see \cite[Appendix]{FE}), provided that $2\sigma >N/2$. (Recall
$(\|u\|_{[1]})^2=\int_0^1|u|^2x^{N/2-1}dx$.) The derivative with respect to $t$
can be treated more simply. Thus we see the equivalence
$$\frac{1}{C}\|u\|_{[1]n}^{(2)}
\leq\|u\|_{[1]n}^{(\infty)}
\leq C\|u\|_{[1], n+s}^{(2)} $$
with $2s>1+N/2$. $\blacksquare$\\

We put
\begin{align*}
\|u\|_n^{(\infty)}&:=\sup_{\mu=0,1}\|u\|_{[\mu]n}^{(\infty)}, \\
\|u\|_n^{(2)}&:=\Big(\sum_{\mu=0,1}(\|u\|_{[\mu]n}^{(2)})\Big)^{1/2}.
\end{align*}

\begin{Proposition}
The gradings $(\|\cdot\|_n^{(\infty)})_n$, $(\|\cdot\|_n^{(2)})_n$ are equivalent and 
the space $\mathfrak{E}$ becomes a tame space by these gradings.
\end{Proposition}

Proof. In fact $\mathfrak{E}$ is a tame direct summand of the
Cartesian product $\mathfrak{E}_{[0]}\times
\mathfrak{E}_{[1]}$, which is a tame space,
(see \cite[p.136, 1.3.1. and 1.3.4]{Hamilton}), by the linear mapping
$L: \mathfrak{E}\rightarrow \mathfrak{E}_{[0]}\times\mathfrak{E}_{[1]}$
$: u\mapsto (u^{[0]}, u^{[1]})$ and
$M:\mathfrak{E}_{[0]}\times\mathfrak{E}_{[1]}\rightarrow\mathfrak{E}:$
$(u_0,u_1)\mapsto u_0+u_1$. Clearly $L$ is tame and $ML=Id_{\mathfrak{E}}$. To verify the tameness of $M$, we note that,
if the support of a function $u$ is included in $[1/6, 5/6]$, then
$$\|\triangle_{[\mu]}^mu\|_{L^{\infty}}
\leq C\sum_{0\leq k\leq m}
\|\triangle_{[1-\mu]}^ku\|_{L^{\infty}}.$$
A proof can be found in \cite[Appendix B]{OJM}. If
$h_{\mu}\in \mathfrak{E}_{[\mu]}$ for $\mu=0,1$, then
$h=M(h_0,h_1)=h_0+h_1$ and
$h^{[0]}=\omega\cdot h_0+\omega\cdot h_1$.
Then by \cite[Proposition 4]{FE} we have
\begin{align*}
\|\triangle_{[0]}^mh^{[0]}\|_{L^{\infty}}&
\leq C\sum_{k\leq m}\|\triangle_{[0]}^kh_0\|_{L^{\infty}}+
\|\triangle_{[0]}^m(\omega\cdot h_1)\|_{L^{\infty}} \\
&\leq C\sum_{k\leq m}\|\triangle_{[0]}^kh_0\|_{L^{\infty}}+
C'\sum_{k\leq m}\|\triangle_{[1]}^k(\omega\cdot h_1)\|_{L^{\infty}} \\
&\leq C\sum_{k\leq m}\|\triangle_{[0]}^kh_0\|_{L^{\infty}}+
C''\sum_{k\leq m}\|\triangle_{[1]}^kh_1\|_{L^{\infty}},
\end{align*}
since the support of $\omega\cdot h_1$ is included in
$[1/6, 2/3]$. Similarly we get an estimate 
 $|\triangle_{[1]}^mh^{[1]}\|_{L^{\infty}}$. These estimates imply the tameness of $M$. Thus 
$\mathfrak{E}$ is a tame direct summand. $\blacksquare$\\

Then $\mathfrak{E}\times\mathfrak{E}$ is a tame space as the
Cartesian product of the tame space $\mathfrak{E}$ and its copy.
Let us denote
\begin{align*}
\|\mathsf{w}\|_n^{(\infty)}&=
\max(\|\tilde{y}\|_n^{(\infty)}, \|\tilde{v}\|_n^{(\infty)}), \\
\|\mathsf{w}\|_n^{(2)}&=
\Big((\|\tilde{y}\|_n^{(2)})^2+(\|\tilde{v}\|_n^{(2)})^2\Big)^{1/2}.
\end{align*}
Then we can claim
\begin{Proposition}
If $2s >1+max(N,5)/2$, then we have
$$\frac{1}{C}\|\mathsf{w}\|_n^{(2)}
\leq\|\mathsf{w}\|_n^{(\infty)}\leq C
\|\mathsf{w}\|_{n+s}^{(2)}.
$$
\end{Proposition}

Since $\mathfrak{P}(\mathsf{w})$ is a smooth function of
$t,x,\tilde{y}^{[0]}, \tilde{y}^{[1]}, D\tilde{y}^{[0]}, D\tilde{y}^{[1]} $,
$\triangle_{[0]}\tilde{y}^{[0]}, \triangle_{[1]}\tilde{y}^{[1]}$,
$\tilde{v}^{[0]}, \tilde{v}^{[1]}, D\tilde{v}^{[0]}, D\tilde{v}^{[1]}$,
$\partial_t\tilde{y}^{[0]}, \partial_t\tilde{y}^{[1]},
\partial_t\tilde{v}^{[0]}, \partial_t\tilde{v}^{[1]}$, thanks to
\cite[Proposition 3, 4, 5]{FE} we see 
$$\|\mathfrak{P}(\mathsf{w})\|_n^{(\infty)}
\leq C(1+\|\mathsf{w}\|_{n+1}^{(\infty)}),$$
that is, $\mathfrak{P}$ is a tame mapping from $\mathfrak{U} \subset
\mathfrak{E}_0\times\mathfrak{E}_0$ into
$\mathfrak{E}\times\mathfrak{E}$.\\

Besides these gradings, we shall use another gradings. Namely, we put
$$\langle u\rangle _{[\mu]\ell}:=
\begin{cases}
\|(\triangle_{[\mu]})^mu\|_{[\mu]} \quad \mbox{for}\quad \ell=2m \\
\|\dot{D}_{[\mu]}(\triangle_{[\mu]})^mu\|_{[\mu]} \quad\mbox{for}\quad
\ell=2m+1,
\end{cases}
$$
where $\mu=0,1$, and
$$\dot{D}_{[0]}=\sqrt{x}\frac{d}{dx},\qquad
\dot{D}_{[1]}=\sqrt{X}\frac{d}{dX} \quad\mbox{with}\quad X=1-x,$$
and for $\mathsf{h}=(h,k)^T$ we put
\begin{align*}
\|\mathsf{h}\|_{[\mu]k}&:=\Big(\sum_{0\leq\ell\leq k}
(\langle h\rangle_{[\mu],\ell+1})^2+
(\langle k\rangle_{[\mu]\ell})^2\Big)^{1/2}, \\
\|\mathsf{h}\|_{[\mu]n}^{\langle\tau\rangle}&:=
\sup_{0\leq t\leq\tau}\sum_{j+k\leq n}\|\partial_t^j\mathsf{h}(t,\cdot)\|_{[\mu]k}, \\
|\|\mathsf{h}\||_{[\mu]n}&:=
\Big(\sum_{j+k\leq n}\int_0^{T}(\|\partial_t^j\mathsf{h}
\|_{[\mu]k})^2dt\Big)^{1/2}, \\
\|\mathsf{h}\|_k&:=\Big((\|\mathsf{h}^{[0]}\|_{[0]k})^2+
(\|\mathsf{h}^{[1]}\|_{[1]k})^2\Big)^{1/2}, \\
\|\mathsf{h}\|_n^{\langle\tau\rangle}&:=
\Big((\|\mathsf{h}^{[0]}\|_{[0]n}^{\langle\tau\rangle})^2+
(\|\mathsf{h}^{[1]}\|_{[1]n}^{\langle\tau\rangle})^2\Big)^{1/2},\\
|\|\mathsf{h}\||_n&:=\Big((|\|\mathsf{h}^{[0]}\||_{[0]n})^2+
(|\|\mathsf{h}^{[1]}\||_{[1]n})^2\Big)^{1/2}.
\end{align*}
Thanks to \cite[Proposition 7]{FE} and the Sobolev's
imbedding with respect to $t$ we can claim
\begin{Proposition}
For any non-negative integer $m$, we have
$$\frac{1}{C}|\|\mathsf{h}\||_n\leq
\|\mathsf{h}\|_n^{\langle T\rangle}
\leq C |\|\mathsf{h}\||_{n+1}
$$
and
$$\frac{1}{C}\|\mathsf{h}\|_m^{(2)}\leq
|\|\mathsf{h}\||_{2m}\leq C
\|\mathsf{h}\|_m^{(2)}.$$
\end{Proposition}

Moreover we put
\begin{align*}
|\mathsf{h}|_n&:=\max_{j+k\leq n, \mu=0,1}
\|\partial_t^j(\dot{D}_{[\mu]})^k\mathsf{h}^{[\mu]}\|_{L^{\infty}}, \\
|\mathsf{h}|_n^{\langle\tau\rangle}&:=\sup_{0\leq t\leq \tau}
|\mathsf{h}(t)|_n.
\end{align*}
By interpolations, we can claim
\begin{Proposition}
For any non-negative integer $m$ we have

$$\frac{1}{C}\|\mathsf{h}\|_m^{(\infty)}\leq
|\mathsf{h}|_{2m}^{\langle T\rangle}
\leq C\|\mathsf{h}\|_m^{(\infty)}.$$
\end{Proposition}

\section{Analysis of the Fr\'{e}chet derivative I}

We have to analyze the Fr\'{e}chet derivative
$D\mathfrak{P}(\mathsf{w})$ of the mapping $\mathfrak{P}$
at a given $\mathsf{w}=(\tilde{y},\tilde{v})^T\in \mathfrak{U}$.
For $\mathsf{h}=(h,k)^T$ we have $D\mathfrak{P}(\mathsf{w})\mathsf{h}
=((DP)_1, (DP)_2)^T$, where
\begin{align*}
(DP)_1&=\frac{\partial h}{\partial t}-Jk-\Big((\partial_yJ)v+(\partial_zJ)vx
\frac{\partial}{\partial x}\Big)h, \\
(DP)_2&=\frac{\partial k}{\partial t}+H_1\mathcal{L}h + \\
&+\Big((\partial_yH_1)\mathcal{L}y+\partial_yH_2+
((\partial_zH_1)\mathcal{L}y+\partial_zH_2)
x\frac{\partial}{\partial x}\Big)h + \\
&+\Big((\partial_vH_1)\mathcal{L}y+\partial_vH_2+
(\partial_wH_2)x\frac{\partial}{\partial x}\Big)k.
\end{align*}

Thanks to the assumption {\bf (B3)} there are analytic functions 
$a_{01}, a_{00}, a_{11}, a_{10}, a_{21}, a_{20}$ of
$x,y(=y^*+\tilde{y}), Dy, D^2y, v(=v^*+\tilde{v}), Dv $, where $D=\partial/\partial x$, on a neighborhood of 
$[0,1]\times \{0\}\times\cdots \times \{0\}$ such that
\begin{subequations}
\begin{align}
(DP)_1&=\frac{\partial h}{\partial t}-Jk+
(a_{01}x(1-x)D+a_{00})h, \\
(DP)_2&=\frac{\partial k}{\partial t}+
H_1\mathcal{L}h+
(a_{11}x(1-x)D+a_{10})h+ \nonumber \\
&+(a_{21}x(1-x)D+a_{20})k.
\end{align}
\end{subequations}

We put
\begin{align*}
\mathfrak{X}&:=L^2([0,1]; x^{3/2}(1-x)^{N/2-1}dx), \\
\mathfrak{X}^1&:=\{\phi\in\mathfrak{X} |\  
\dot{D}\phi:=\sqrt{x(1-x)}\frac{d\phi}{dx}\in\mathfrak{X}\}, \\
\mathfrak{X}^2&:=\{\phi\in \mathfrak{X}^1 |\  -\Lambda\phi\in\mathfrak{X}\},
\end{align*}
with
\begin{equation}
\Lambda=x(1-x)\frac{d^2}{dx^2}+\Big(\frac{5}{2}(1-x)-
\frac{N}{2}\Big)\frac{d}{dx},
\end{equation}

We claim
\begin{Proposition}
Given $\mathsf{g}\in C([0,T], \mathfrak{X}^1\times
\mathfrak{X})$, the
equation $D\mathfrak{P}(\mathsf{w})\mathsf{h}=\mathsf{g}$
admits a unique solution $\mathsf{h}\in C([0,T], \mathfrak{X}^2\times\mathfrak{X}^1)$ such that
$\mathsf{h}(t=0)=(0,0)^T$.
\end{Proposition}

Proof. We can rewrite
$$(DP)_2=\frac{\partial k}{\partial t}-H_1\Lambda h
+b_1\check{D}h+b_0h+a_{21}\check{D}k+a_{20}k,
$$
where
\begin{align*}
\check{D}&:=x(1-x)D, \\
b_1&:=H_1\ell_1+a_{11}, \\
b_0&:=H_1L_0+a_{10}.
\end{align*}
Of course $b_1, b_0$ are also analytic on a neighborhood of
$[0,1]\times\{0\}\times\cdots$.
Then we can write the equation
$D\mathfrak{P}(\mathsf{w})\mathsf{h}=\mathsf{g}(=(g_1,g_2)^T)$ as
\begin{equation}
\frac{\partial}{\partial t}
\begin{bmatrix}
h \\
k
\end{bmatrix}
+
\begin{bmatrix}
\mathfrak{a}_1 & -J \\
\mathcal{A} & \mathfrak{a}_2
\end{bmatrix}
\begin{bmatrix}
h \\
k
\end{bmatrix}
=
\begin{bmatrix}
g_1 \\
g_2
\end{bmatrix}.
\label{Eq16}
\end{equation}
Here
\begin{align*}
\mathfrak{a}_1&=a_{01}\check{D}+a_{00}, \\
\mathfrak{a}_2&:=a_{21}\check{D}+a_{20}, \\
\mathcal{A}&:=-H_1\Lambda + b_1\check{D}+b_0.
\end{align*}
The standard calculation gives
\begin{align}
&\frac{1}{2}\frac{d}{dt}\Big(\|k\|^2+((H_1/J)\dot{D}h|\dot{D}h)\Big) +\nonumber\\
&+(\beta_1\dot{D}h|\dot{D}h)+
(\beta_2\dot{D}h|h)+(\beta_3\dot{D}h|k)+
+(\beta_4h|k)+(\beta_5k|k)= \nonumber\\
&=((H_1/J)\dot{D}h|\dot{D}g_1)+(k|g_2), \label{Eneqy}
\end{align}
where $\dot{D}=\sqrt{x(1-x)}D$ and
\begin{align*}
\beta_1&=-\frac{1}{4}(3+(N+3)x+2\check{D})(H_1/J)a_{01}-\frac{1}{2}
\frac{\partial(H_1/J)}{\partial t} +
(H_1/J)(\check{D}a_{01}+a_{00}), \\
\beta_2&=(H_1/J)\dot{D}a_{00}, \\
\beta_3&=-(H_1/J)\dot{D}J+\dot{D}H_1+
\sqrt{x(1-x)}(b_1+a_{21}), \\
\beta_4&=b_0,\qquad \beta_5=a_{20}.
\end{align*}
Of course $(\cdot|\cdot)$ and $\|\cdot\|$ stand for the inner product and the norm of the Hilbert space $\mathfrak{X}$. We have used the following
formulas:\\

{\bf Formula 1}: {\it If $\phi\in\mathfrak{X}^2, \psi\in\mathfrak{X}^1,
\alpha\in C^1([0,1])$, then
\begin{equation}
(-\alpha\Lambda\phi|\psi)=
(\alpha\dot{D}\phi|\dot{D}\psi)+
((D\alpha)\check{D}\phi|\psi). \label{F1}
\end{equation}}\\

{\bf Formula 2}: {\it If $\phi \in\mathfrak{X}^2$ and $\alpha\in
C^1([0,1])$, then
\begin{equation}
(\alpha\dot{D}\phi|\dot{D}\check{D}\phi)=
(\alpha^*\dot{D}\phi|\dot{D}\phi), 
\label{F2}
\end{equation}
with
$$\alpha^*=-\frac{1}{4}(3+(N+3)x+2\check{D})\alpha.$$ }\\

Since $\mathsf{w}$ is confined to $\mathfrak{U}$,
we can assume
$$\frac{1}{M_0}<J<M_0,\qquad \frac{1}{M_0}<H_1<M_0 $$
with a constant $M_0$ independent of  $\mathsf{w}$ thanks to the assumption {\bf (B2)}. Now the energy
$$\mathcal{E}:=\|k\|^2+((H_1/J)\dot{D}h|\dot{D}h) $$
enjoys the inequality
$$\frac{1}{2}\frac{d\mathcal{E}}{dt}\leq
M(\|\mathsf{h}\|_{\mathfrak{H}}^2+
\|h\|_{\mathfrak{H}}\|\mathsf{g}\|_{\mathfrak{H}}),
$$
where $\mathfrak{H}=\mathfrak{X}^1\times\mathfrak{X}$ and
$$
\|(\phi,\psi)^T\|_{\mathfrak{H}}^2=\|\phi\|_{\mathfrak{X}^1}^2+
\|\psi\|_{\mathfrak{X}}^2
=\|\phi\|^2+\|\dot{D}\phi\|^2+\|\psi\|^2$$
and
$$M=\sum_{j=1}^5\|\beta_j\|_{L^{\infty}}+(M_0)^2+1.
$$
Since $\mathcal{E}$ is equivalent to $\|k\|^2+\|\dot{D}h\|^2$, the Gronwall
argument and application of the Kato's theory (\cite{Kato11})
deduce the conclusion. Here $\|h\|$ should be estimated by $\mathcal{E}$
as follows: The first component of (\ref{Eq16}) implies
$$h(t)=\int_0^t(-a_{01}\check{D}h-a_{00}h+Jk+g_1)(t')dt', $$
therefore
$$\|h(t)\|\leq C\int_0^t\|h(t')\|dt'+
\int_0^t(C\mathcal{E}(t')+\|g_1(t')\|)dt',$$
where
$C=\max(\|a_{00}\|_{L^{\infty}},\|a_{01}\|_{L^{\infty}}M_0^2+ M_0)$, which implies, through the Gronwall's argument,
$$\|h(t)\|\leq\int_0^t
(e^{C(t-t')}-1)(C\mathcal{E}(t')+\|g_1(t')\|)dt'.$$
As the result the solution enjoys
$$\|\mathsf{h}(t)\|_{\mathfrak{H}}
\leq C\int_0^t\|\mathsf{g}(t')\|_{\mathfrak{H}}dt'.
$$
$\blacksquare$\\

Here, in order to make sure, let us sketch proofs of {\bf Formula 1},
{\bf Formula 2}, and  (\ref{Eneqy}).\\

Proof of (\ref{F1}): if $\psi\in \mathfrak{X}^1$, then 
$$\psi(1/2)+\int_{1/2}^x\frac{\dot{D}\psi(x')}{\sqrt{x'(1-x')}}dx'$$
implies
$$|\psi(x)|\leq C x^{-3/4}(1-x)^{-N/4+1/2},$$
and, if $\phi\in\mathfrak{X}^2$, then
\begin{align*}
x^{5/2}(1-x)^{N/2}\frac{d\phi}{dx}&=
x^{5/2}(1-x)^{N/2}\frac{d\phi}{dx}\Big|_{x=1/2} +\\
&-\int_{1/2}^x
\Lambda\phi(x')x'^{3/2}(1-x')^{N/2-1}dx'
\end{align*}
implies
$$
\Big|\frac{d\phi}{dx}\Big|\leq C
x^{-5/4}(1-x)^{-N/4}.$$
Actually the finite constant
$$x^{5/2}(1-x)^{N/2}D\phi|_{x=1/2}+
\int_0^{1/2}\Lambda\phi(x')x'^{3/2}(1-x')^{N/2-1}dx'$$
should vanish in order to $\dot{D}\phi \in\mathfrak{X}$ and so on. 
Therefore the boundary terms in the integration by parts vanish
at $x=+0, 1-0$ and we get the desired equality.\\

Proof of (\ref{F2}): We see
\begin{align*}
(\alpha\dot{D}\phi|\dot{D}\check{D}\phi)&=
\int_0^1\alpha x(1-x)(D\phi)D(x(1-x)D\phi)x^{3/2}(1-x)^{N/2-1}dx \\
&=I +(\alpha(1-2x)\dot{D}\phi|\dot{D}\phi),
\end{align*}
where
\begin{align*}
I&:=\int_0^1\alpha(x(1-x))^2(D\phi)(D^2\phi)x^{3/2}(1-x){N/2-1}dx \\
&=\int_0^1\frac{\alpha}{2}D(D\phi)^2x^{7/2}(1-x)^{N/2+1}dx \\
&=-\int_0^1
D\Big(\frac{\alpha}{2}x^{7/2}(1-x)^{N/2+1}\Big)(D\phi)^2dx.
\end{align*}
Here the integration by parts has been done by using
$$|D\phi|\leq C x^{-5/4}(1-x)^{-N/4}$$
which holds for $\phi\in\mathfrak{X}^2$. 
Then we see
$$I=-\Big(\Big(\frac{\check{D}\alpha}{2}+
\frac{\alpha}{2}\Big(\frac{7}{2}(1-x)+\Big(\frac{N}{2}+1\Big)x\Big)\Big)
\dot{D}\phi\Big|\dot{D}\phi\Big),$$
and get (\ref{F2}). \\

Proof of (\ref{Eneqy}): multiplying the second component of the equation
(\ref{Eq16}) by $k$ and integrating it, we get
\begin{align*}
\frac{1}{2}\frac{d}{dt}
\|k\|^2-(H_1\Lambda h|k)&+(b_1\check{D}h|k)+(b_0h|k) + \\
&+(a_{21}\check{D}h|k)+(a_{20}k|k)=(g_2|k).
\end{align*}
By {\bf Formula 1} we see
\begin{align}
\frac{1}{2}\frac{d}{dt}
\|k\|^2
+(H_1\dot{D}h|\dot{D}k)&
+((DH_1)\check{D}h|k)
+(b_1\check{D}h|k)+(b_0h|k) + \nonumber\\
&+(a_{21}\check{D}h|k)+(a_{20}k|k)=(g_2|k). \label{*}
\end{align}
On the other hand, operating $\dot{D}$ on the first component of
(\ref{Eq16}), we get
\begin{align*}
\dot{D}k&=\frac{1}{J}\partial_t(\dot{D}h)+
\frac{a_{01}}{J}\dot{D}\check{D}h+ \\
&+\frac{1}{J}(\check{D}a_{01}+a_{00})\dot{D}h+
\frac{\dot{D}a_{00}}{J}h-\frac{\dot{D}J}{J}h-\frac{\dot{D}J}{J}k-
\frac{1}{J}\dot{D}g_1.
\end{align*}
Inserting this into the second term of the left-hand
side of (\ref{*}), we get
\begin{align*}
&\frac{1}{2}\frac{d}{dt}\|k\|^2+\frac{1}{2}\frac{d}{dt}
((H_1/J)\dot{D}h|\dot{D}h)+
((H_1a_{01}/J)\dot{D}h|\dot{D}\check{D}h)+ \\
&+(((1/2)(H_1/J)_t+(H_1/J)(\check{D}a_{01}+a_{00})\dot{D}h|\dot{D}h)+
((H_1/J)(\dot{D}a_{00})(\dot{D}h)|h) + \\
&-(H_1\dot{D}J/J)\dot{D}h|k)-
((H_1/J)\dot{D}h|\dot{D}g_1) + \\
&+((b_1+a_{21})\check{D}h|k)+(b_0h|k)+(a_{20}k|k)=(g_2|k).
\end{align*}
Applying {\bf Formula 2} to the third term of the left-hand
side, we get the desired (\ref{Eneqy}).

\section{Analysis of the Fr\'echet derivative, II}

We are going to show the smoothness of the solution $\mathsf{h}$
of the equation $D\mathfrak{P}(\mathsf{w})\mathsf{h}=\mathsf{g}$ 
for $\mathsf{g}\in \mathfrak{E}\times\mathfrak{E}$ given and to get a tame estimate of the mapping $(\mathsf{w}, \mathsf{g}) \mapsto \mathsf{h}$.
Suppose that $\mathsf{w}, \mathsf{g} \in C^{\infty}([0,T]\times[0,1])$ and let $\mathsf{h}\in C([0,1], \mathfrak{X}^2\times \mathfrak{X}^1)$ be
he solution of $D\mathfrak{P}(\mathsf{w})\mathsf{h}=\mathsf{g},
\mathsf{h}|_{t=0}=\mathsf{0}$.\\

We are going to show $\mathsf{h}\in C^{\infty}([0,T]\times [0,1])$. Using
the cut off function $\omega\in C^{\infty}(\mathbb{R})$ such that
$\omega(x)=1$ for $x\leq 1/3$, $0<\omega<1$ for
$1/3<x<2/3$ and $\omega(x)=0$ for $2/3\leq x$, we put
$$\mathsf{h}^{[0]}(t,x)=
\begin{bmatrix}
h^{[0]} \\
k^{[0]}
\end{bmatrix}
=\omega(x)\mathsf{h}(t,x), \qquad
\mathsf{h}^{[1]}(t,x)=
\begin{bmatrix}
h^{[1]} \\
k^{[1]}
\end{bmatrix}
=(1-\omega(x))\mathsf{h}(t, x),
$$
and $$\mathsf{g}^{[0]}(t,x)=
\begin{bmatrix}
g_1^{[0]} \\
g_2^{[0]}
\end{bmatrix}
=\omega\cdot\mathsf{g}(t,x)$$ and so on. Then the equation (\ref{Eq16}) turns out to be
\begin{equation}
\frac{\partial}{\partial t}
\begin{bmatrix}
h^{[\mu]} \\
k^{[\mu]}
\end{bmatrix}+
\begin{bmatrix}
\mathfrak{a}_1^{[\mu]} & -J \\
\mathcal{A}^{[\mu]} & \mathfrak{a}_2^{[\mu]}
\end{bmatrix}
\begin{bmatrix}
h^{[\mu]} \\
k^{[\mu]}
\end{bmatrix}
=\begin{bmatrix}
g_1^{\mu]} \\
g_2^{\mu]}
\end{bmatrix}
+(-1)^{\mu}
\begin{bmatrix}
c_{11} & 0 \\
\mathfrak{c}_{21} & c_{22}
\end{bmatrix}
\begin{bmatrix}
h^{[1-\mu]} \\
k^{[1-\mu]}
\end{bmatrix},
\quad \mu=0,1,
\label{Eq301}
\end{equation}
where
\begin{align*}
\mathfrak{a}_1^{[\mu]}&=a_{01}\check{D}+a_{00}
-(-1)^{\mu}a_{01}\check{D}\omega, \\
\mathfrak{a}_2^{[\mu]}&=a_{21}\check{D}+a_{20}
-(-1)^{\mu}a_{21}\check{D}\omega, \\
\mathcal{A}^{[\mu]}&=-H_1\Lambda+
(b_1+(-1)^{\mu}2H_1(D\omega))\check{D}+b_0+
(-1)^{\mu}(H_1\Lambda-b_1\check{D})\omega, \\
c_{11}&=a_{01}\check{D}\omega, \\
\mathfrak{c}_{21}&=-2H_1(D\omega)\check{D}+
b_1\check{D}\omega-H_1\Lambda\omega, \\
c_{22}&=a_{21}\check{D}\omega.
\end{align*}
We can write
$$\mathcal{A}^{[\mu]}=-b_2^{[\mu]}\triangle_{[\mu]}+b_1^{[\mu]}\check{D}_{[\mu]}+b_0^{[\mu]},
$$
where
\begin{align*}
&\check{D}_{[0]}=x\frac{\partial}{\partial x},
\qquad
\check{D}_{[1]}=X\frac{\partial}{\partial X}\quad\mbox{with}\quad X=1-x, \\
& b_2^{[0]}=H_1\cdot (1-x), \qquad b_1^{[0]}=\frac{N}{2}H_1+
(b_1+2H_1(D\omega))(1-x), \\
&b_2^{[1]}=H_1\cdot x, \qquad b_1^{[1]}=\frac{5}{2}H_1-
(b_1-2H_1(D\omega))x, \\
&b_0^{[\mu]}=b_0-(-1)^{\mu}(H_1\Lambda -b_1\check{D})\omega.
\end{align*}
 We have $1/C\leq b_2^{[\mu]} < C$ on
$I_{[\mu]}$, where $I_{[0]}=[0,2/3], I_{[1]}=[1/3,1]$. We note that the 
regularity of $\mathsf{h}$ can be reduced to that of
$\mathsf{h}^{[0]}, \mathsf{h}^{[1]}$. In fact, if
$\mathsf{h}^{[0]} \in C^{\infty}([0,T]\times[0, 2/3])$, then
$\mathsf{h}(t,x)=\mathsf{h}^{[0]}(t,x)/\omega(x)$ is smooth on 
$0\leq t\leq T, 0\leq x <2/3$ and so on. But
the regularity of the solution of (\ref{Eq301}) can be proved
by Kato's theory in \cite{OJM12}. Namely we put
\begin{align*}
&\hat{\mathfrak{H}}=\mathfrak{H}_{[0]}\times
\mathfrak{H}_{[1]}\times\mathbb{R}, \\
&\mathfrak{H}_{[\mu]}=\mathfrak{X}_{[\mu]0}^1
\times\mathfrak{X}_{[\mu]},\quad \mu=0,1, 
\end{align*}
\begin{align*}
&\mathfrak{X}_{[0]}=L^2([0,2/3], x^{3/2}dx), \\
&\mathfrak{X}_{[0]}^1=\{\phi\in\mathfrak{X}_{[0]} |\  
\dot{D}_{[0]}\phi:=\sqrt{x}\frac{d\phi}{dx}\in\mathfrak{X}_{[0]}\},
\quad \mathfrak{X}_{[0]0}^1=\{\phi\in\mathfrak{X}_{[0]} |\  \phi(2/3)=0 \}, \\
&\mathfrak{X}_{[0]}^2=\{\phi\in\mathfrak{X}_{[0]}
|\  \triangle_{[0]}\phi\in\mathfrak{X}_{[0]}\},
\quad
\mathfrak{X}_{[0](0)}^2=\mathfrak{X}_{[0]}^2\cap\mathfrak{X}_{[0](0)}^1, 
\end{align*}
\begin{align*}
&\mathfrak{X}_{[1]}=L^2([1/3,1],(1-x)^{N/2-1}dx ), \\
&\mathfrak{X}_{[1]}^1=\{\phi\in\mathfrak{X}_{[1]} |\  \dot{D}_{[1]}\phi:=
-\sqrt{1-x}\frac{d\phi}{dx}\in\mathfrak{X}_{[1]} \}, \quad \mathfrak{X}_{[1]0}^1=\{\phi\in\mathfrak{X}_{[1]}^1 |\  \phi(1/3)=0 \}, \\
&\mathfrak{X}_{[1]}^2=\{\phi\in\mathfrak{X}_{[1]}^1 |\  
\triangle_{[1]}\phi \in\mathfrak{X}_{[1]} \},
\quad \mathfrak{X}_{[1](0)}^2=
\mathfrak{X}_{[1]}^2\cap\mathfrak{X}_{[1]0}^1, 
\end{align*}
\begin{align*}
&D(\hat{\mathfrak{A}}(t))=\hat{\mathfrak{G}}=
\mathfrak{G}_{[0]}\times\mathfrak{G}_{[1]}\times\mathbb{R}, \\
&\mathfrak{G}_{[\mu]}=\mathfrak{X}_{[\mu](0)}^2\times
\mathfrak{X}_{[\mu]0}^1,\quad \mu=0,1, \\
&\hat{\mathfrak{A}}(t)=\mathfrak{A}_{[0]}(t)\otimes
\mathfrak{A}_{[1]}(t)\otimes 0+B(t), \\
&\mathfrak{A}_{[\mu]}(t)=
\begin{bmatrix}
\mathfrak{a}_1^{[\mu]} & -J \\
\mathcal{A}^{[\mu]} & \mathfrak{a}_2^{[\mu]}
\end{bmatrix},
\quad \mu=0,1, 
\end{align*}
$$
B(t)=\begin{bmatrix}
0 & 0 & -c_{11} & 0& -g_1^{[0]} \\
0 & 0 & -\mathfrak{c}_{21} & -c_{22} & -g_2^{[0]} \\
c_{11} & 0 & 0&0&-g_1^{[1]} \\
\mathfrak{c}_{21} & c_{22} & 0 & 0 & -g_2^{[1]} \\
0&0&0&0&0 
\end{bmatrix}.
$$
Then the equation (\ref{Eq301}) is reduced to
$$\frac{\partial \mathsf{U}}{\partial t}+\hat{\mathfrak{A}}(t)\mathsf{U}=\mathsf{0} $$
for $\mathsf{U}=(h^{[0]}, k^{[0]}, h^{[1]}, k^{[1]}, 1)^T$ and the initial condition reads
$$\mathsf{U}|_{t=0}=\phi_0=(0,0,0,0,1)^T.$$
The same application of the Kato's theory as 
\cite[pp.569-571]{OJM} guarantees 
$\mathsf{h}^{[0]}\in C^{\infty}([0,T]\times[0,2/3]),
\mathsf{h}^{[1]}\in C^{\infty}
([0,T]\times[1/3,1])$ provided that
$\mathsf{g}\in C^{\infty}([0,T]\times[0,1])$.

Here the ellipticity of $\mathfrak{A}(t)=\mathfrak{A}_{[\mu]}(t),
\mu=0,1$, reads
\begin{equation}
\|\mathsf{h}\|_{n+1}\leq C(\|\mathfrak{A}\mathsf{h}\|_n+
(1+|\mathsf{a}|_{n+4})\|\mathsf{h}\|_0),
\label{Eq302}
\end{equation}
with
\begin{align*}
\|\mathsf{h}\|_k&=\|\mathsf{h}\|_{[\mu]k}:=
\Big(\sum_{0\leq\ell\leq k}\langle h\rangle_{[\mu],\ell+1}^2+
\langle k\rangle_{[\mu]\ell}^2\Big)^{1/2}, \\
\langle\phi\rangle_{[\mu]\ell}&=
\begin{cases}
\|\triangle_{[\mu]}^m\phi\|_{\mathfrak{X}_{[\mu]}} &\quad\mbox{for}\quad \ell=2m \\
\|\dot{D}_{[\mu]}\triangle_{[\mu]}^m\phi\|_{\mathfrak{X}_{[\mu]}} &\quad
\mbox{for}\quad \ell=2m+1
\end{cases},
\\
\mathsf{a}&=\mathsf{a}_{[\mu]}=(a_i^{[\mu]})_{i=0}^7=
(b_0^{[\mu]},b_1^{[\mu]}, b_2^{[\mu]}, a_{01}, a_{00}, a_{21},a_{20},J), \\
|\mathsf{a}|_n&=\max_{j+k\leq n, 0\leq i\leq 7}
\|\partial_t^j\dot{D_{[\mu]}}^ka_i^{[\mu]}\|_{L^{\infty}}.
\end{align*}
In fact, for $\mathfrak{a}_1=\mathfrak{a}_1^{[\mu]},
\mathfrak{a}_2=\mathfrak{a}_2^{[\mu]}, \mu=0,1$, we have
\begin{align*}
\|\mathfrak{a}_1h\|_1&\leq C(\varepsilon\|h\|_2+\|h\|_1), \\
\|\mathfrak{a}_2k\|_0&\leq C(\varepsilon\|k\|_1+\|k\|_0),
\end{align*}
since $a_{01}=\partial_zJv/(1-x),
a_{21}=\partial_wH_2/(1-x)$ enjoy $|a_{01}|+|a_{21}|\leq \varepsilon(\ll 1)$
provided that the domain $U$ of $y, z, v, w$ is sufficiently small. 
On the other hand the estimates of the commutators for
$\mathcal{A}=\mathcal{A}^{[\mu]}, \mathfrak{a}=\mathfrak{a}_1^{[\mu]},
\mathfrak{a}_2^{[\mu]}, J$ are
\begin{align*}
\|[\triangle,\mathcal{A}]\phi\|_n &\leq
C(|\mathsf{a}|_2\|\phi\|_{n+3}+|\mathsf{a}|_{n+5}\|\phi\|_0), \\
\|[\triangle,\mathfrak{a}]\phi\|_n &\leq
C(|\mathsf{a}|_3\|\phi\|_{n+2}+|\mathsf{a}|_{n+5}\|\phi\|_0), \\
\|[\triangle,J]\phi\|_n &\leq
C(|\mathsf{a}|_4\|\phi\|_{n+1}+|\mathsf{a}|_{n+5}\|\phi\|_0).
\end{align*}
Thus we can derive the elliptic estimate for $\mathfrak{A}=\mathfrak{A}^{[\mu]}, \mu=0,1$.
Then the ellipticity of $\hat{\mathfrak{A}}(t)$ can be verified since the supports of the 
coefficients of the interacting terms
$c_{11}, \mathfrak{c}_{21}, c_{22}$ are included in the interval $[1/3,2/3]$
away from the both singular boundary points.

The compatibility can be verified using the fact that
$\mathsf{g}^{[0]}(t,x)$ and all its higher derivatives with respect to
$x$ vanish at $x=2/3$ and $\mathsf{g}^{[1]}$ and its higher derivatives vanish, too, at
$x=1/3$. (See \cite[pp. 570-571]{OJM}.)\\

Now, we are going to get a tame estimate of the mapping
$(\mathsf{w},\mathsf{g}) \mapsto \mathsf{h}$. 

First we consider the component $\mathsf{h}=\mathsf{h}^{[\mu]}, \mu=0,1$, which satisfies
$$\frac{\partial\mathsf{h}}{\partial t}+
\mathfrak{A}\mathsf{h}=\mathsf{f},$$
with
\begin{align*}
&\mathfrak{A}=\mathfrak{A}_{[\mu]}=
\begin{bmatrix}
\mathfrak{a}_1^{[\mu]} & -J \\
\mathcal{A}^{[\mu]} & \mathfrak{a}_2^{[\mu]}
\end{bmatrix}, \\
&\mathfrak{a}_1^{[\mu]}=a_{01}\check{D}+a_{00}^{[\mu]}=a_{01}\check{D}+
a_{00}-(-1)^{\mu}a_{01}\check{D}\omega, \\
&\mathfrak{a}_2^{[\mu]}=a_{21}\check{D}+a_{20}^{[\mu]}=
a_{21}\check{D}+a_{20}-(-1)^{\mu}a_{21}\check{D}\omega, \\
&\mathcal{A}^{[\mu]}=-b_2^{[\mu]}\triangle_{[\mu]}+b_1^{[\mu]}\check{D}_{[\mu]}+b_0^{[\mu]}, \\
&\mathsf{f}=\mathsf{f}^{[\mu]}=
\begin{bmatrix}
f_1^{\mu]} \\
f_2^{[\mu]}
\end{bmatrix}
= \\
&=\begin{bmatrix}
g_1^{[\mu]}+(-1)^{\mu}c_{11}h^{[1-\mu]} \\
g_2^{[\mu]}+(-1)^{\mu}
(c_{21}\check{D}+c_{20})h^{[1-\mu]}+c_{22}k^{[1-\mu]}.
\end{bmatrix}.
\end{align*}

The elliptic estimate reads
$$\|\mathsf{h}\|_{n+1}\leq
C(\|\mathfrak{A}\mathsf{h}\|_n+(1+|\mathsf{a}|_{n+4})\|h\|_0).$$

The energy estimate reads
$$\|\mathsf{H}\|\leq C
(\|\mathsf{H}|_{t=0}\|+\int_0^t\|\mathsf{F}(t')\|dt'), $$
where
$$\|\mathsf{H}\|=(\|{H}\|^2+\|\dot{D}_{[\mu]}{H}\|^2+\|K\|^2)^{1/2}
$$
with $\|\cdot\|=\|\cdot\|_{[\mu]}, \mathsf{H}=(H,K)^T$, for
any solution $\mathsf{H}(t)$ of
$$\frac{\partial\mathsf{H}}{\partial t}+\mathfrak{A}(t)\mathsf{H}=\mathsf{F}(t). $$

Let us put
\begin{align*}
&X(j,k;\mathsf{H})=\|\partial_t^j\mathsf{H}\|_k, \\
&Z(n;\mathsf{H})=\sum_{j+k=n}\|\partial_t^j\mathsf{H}\|_k=\sum_{j+k=n}X(j,k;\mathsf{H}),\\
&W(n;\mathsf{H})=\sum_{j+k\leq n}\|\partial_t^j\mathsf{H}\|_k=
\sum_{j+k\leq n}X(j,k;\mathsf{H})=
\sum_{\nu\leq n}Z(\nu; \mathsf{H}).
\end{align*}
Since $Z(n+1;\mathsf{h})=Z(n;\partial_t\mathsf{h})+
\|\mathsf{h}\|_{n+1}$, an induction on
$n$ leads us to
\begin{equation}
Z(n+1;\mathsf{h})\leq C(
\|\partial_t^{n+1}\mathsf{h}\|+
\sum_{j+k=n}\|\mathsf{f}_j\|_k+
\sum_{j+k=n}
(1+|\mathsf{a}|_{k+4})\|\partial_t^j\mathsf{h}\|),
\label{TE01}
\end{equation}
where
$$\mathsf{f}_j:=\partial_t^j\mathsf{f}-
[\partial_t^j,\mathfrak{A}]\mathsf{h}.$$
Using
\begin{equation}
\|\mathfrak{A}\mathsf{h}\|_k\leq C(\|\mathsf{h}\|_{k+1}+
|\mathsf{a}|_{k+4}\|\mathsf{h}\|)
\label{TE02}
\end{equation}
and interpolations, we have
\begin{equation}
Z(n+1;\mathsf{h})\leq C(
\|\partial_t^{n+1}\mathsf{h}\|+
Z(n;\mathsf{f})+Z(n;\mathsf{h})+|\mathsf{a}|_{n+4}),
\label{TE03}
\end{equation}
provided that $|\mathsf{a}|_5$ and 
$\sum_{\mu=0,1}\|\mathsf{h}\|^{[\mu]}\leq C\sum_{\mu=0,1}
\int_0^t
\|\mathsf{g}^{[\mu]}\|$ are bounded.
The energy estimate for $\mathsf{H}=\partial_t^{n+1}\mathsf{h}$, which
satisfies
$$\partial_t\mathsf{H}+\mathfrak{A}\mathsf{H}=\mathsf{f}_{n+1}
(=\partial_t^{n+1}\mathsf{f}-
[\partial_t^{n+1},\mathfrak{A}]\mathsf{h}),$$
gives
\begin{align}
\|\partial_t^{n+1}\mathsf{h}\|&\leq C(\|\partial_t^{n+1}\mathsf{h}|_{t=0}\|+
\int_0^t\|\partial_t^{n+1}\mathsf{f}\|+ \nonumber \\
&+\int_0^tZ(n+1;\mathsf{h})+|\mathsf{a}|_{n+5}).
\label{TE04}
\end{align}
Substituting (\ref{TE04}) into (\ref{TE03}), we get
$$Z(n+1;\mathsf{h})(t)\leq C\Big(
\int_0^tW(n+1;\mathsf{h})+\mathsf{F}_n),$$
where
$$
\mathsf{F}_n=\|\partial_t^{n+1}\mathsf{h}|_{t=0}\|+
\int_0^t\|\partial_t^{n+1}\mathsf{f}\|+|\mathsf{a}|_{n+5}+Z(n;\mathsf{f})+
Z(n;\mathsf{h}).
$$
A Gronwall's argument implies
\begin{equation}
W(n+1;\mathsf{h})\leq C
(W_0(n+1;\mathsf{h})+\|\mathsf{f}\|_n^{\langle t\rangle}+
\int_0^tW(n+1;\mathsf{f})+|\mathsf{a}|_{n+5}),
\label{TE05}
\end{equation}
where we put
\begin{align*}
&X_0(j,k;\mathsf{H})=\|\partial_t^j\mathsf{H}|_{t=0}\|_k, \\
&W_0(n;\mathsf{H})=\sum_{j+k\leq n}X_0(j,k; \mathsf{H}).
\end{align*}
Here we have eliminated $W(n,\mathsf{h})$ in the right-hand side by
induction on $n$. 

We have to estimate $W_0(n+1;\mathsf{h})$.
Using (\ref{TE02}) 
and interpolations, we can verify
\begin{align*}
&X_0(j+1,k)\leq C(X_0(j,k+1)+
W_0({j+k})+|\mathsf{a}|_{k+4}X_0(j,0)+\\
&+|\mathsf{a}|_{k+j+3}W_0(0;\mathsf{f})+W_0({j+k};\mathsf{f}),
\end{align*}
where $X_0(j,k)=X_0(j,k;\mathsf{h}), W_0(n)=W_0(n;\mathsf{h})$.
Then, by induction on $n$ and $j$, we get the inequality
$$X_0(j,k)\leq C(W_0(n)+|\mathsf{a}|_{n+4}+W_0(n;\mathsf{f}))\quad\mbox{for}\quad j+k=n+1,$$
provided that $W_0(0;\mathsf{f}))=W_0(0;\mathsf{g})$ and
$|\mathsf{a}|_4$ are bounded. This implies
$$W_0({n+1};\mathsf{h})\leq C(W_0(n;\mathsf{h})+
|\mathsf{a}|_{n+4}+W_0(n;\mathsf{f})),$$
and eliminating $W_0(n;\mathsf{h})$ on the right-hand side by
induction on $n$, we get
$$W_0({n+1};\mathsf{h})\leq C(1+|\mathsf{a}|_{n+4}+W_0(n;\mathsf{f})),$$
that is,
$$W_0({n+1};\mathsf{h^{[\mu]}})\leq C(1+
|\mathsf{a}_{[\mu]}|_{n+4}+W_0(n;\mathsf{f^{[\mu]}})),\quad
\mu=0,1.$$
Then, using 
$$W_0(n;\mathsf{f}^{[\mu]})\leq W_0(n;\mathsf{g}^{[\mu]})+
CW_0(n;\mathsf{h}^{[1-\mu]}),$$
the induction on $n$ gives
\begin{equation}
W_0({n+1};\mathsf{h}^{[\mu]})\leq
C(1+W_0(n;\mathsf{g}^{[\mu]})+a_n),\quad \mu=0,1
\label{TE06}
\end{equation}
provided that $|\mathsf{a}_{[\mu]}|_4$ and 
$W_0(0;(\mathsf{f}^{[\mu]})=W_0(0;\mathsf{g}^{[\mu]})$
are bounded,
where $$a_n:=\max\{ |\mathsf{a}_{[\mu]}|_{n+5}^{\langle T\rangle} : \mu=0,1\}.$$ 
Since
$$W_0(n;\mathsf{g}^{[\mu]})\leq C\|\mathsf{g}\|_{n}^{\langle T\rangle},$$
we have
\begin{equation}
W_0(n+1;\mathsf{h}^{[\mu]})\leq
C(1+\|\mathsf{g}\|_{n}^{\langle T\rangle}+a_n).
\label{TE07}
\end{equation}

Substituting (\ref{TE07}) into (\ref{TE05}),
we get the simultaneous inequalities
\begin{align*}
W_{[\mu]}(n+1;\mathsf{h}^{[\mu]})(t)\leq & C
\Big(1+
\int_0^t W_{[1-\mu]}(n+1;\mathsf{h}^{[1-\mu]})(t')dt'+ \\
&+\int_0^t
W_{[\mu]}(n+1;\mathsf{g}^{[\mu]})(t')dt'+
\|\mathsf{g}\|_{n}^{\langle T\rangle}+a_n\Big),\quad \mu=0,1.
\end{align*}
Then the Gronwall's argument 
applied to the quantity
$$Q(t)=
W_{[0]}(n+1,\mathsf{h}^{[0]})+W_{[1]}(n+1;\mathsf{h}^{[1]})$$
gives
\begin{equation}
\|\mathsf{h}\|_{n+1}^{\langle T\rangle}\leq
C(1+\|\mathsf{g}\|_{n}^{\langle T\rangle}+a_n).
\end{equation}
Since $a_n\leq C|\mathsf{w}|_{n+7}^{\langle T\rangle}
$
, we have the desired tame estimate
$$\|\mathsf{h}\|_m^{(2)}\leq C
(1+\|\mathsf{g}\|_m^{(2)}+
\|\mathsf{w}\|_{m+r}^{(2)}),
$$
provided that $2r>7+\max(N,5)/2$.

\section{Corollaries}

Let us describe the spectral properties of the linear
operator $\mathcal{L}$. We write (\ref{2}) as
\begin{equation}
\mathcal{L}=-\frac{1}{b(x)}\frac{d}{dx}a(x)\frac{d}{dx}+L_0(x), 
\end{equation}
where
\begin{subequations}
\begin{align}
a(x)&=x^{5/2}(1-x)^{N/2}M(x), \\
b(x)&=x^{3/2}(1-x)^{N/2-1}M(x), \\
M(x)&=\exp\Big[\int_0^x\ell_1(x')dx'\Big]
\end{align}
\end{subequations}
Note that $M$ is positive and analytic on a neighborhood of
$[0,1]$. The Liouville transformation
\begin{subequations}
\begin{align}
\xi&=\int_{1/2}^x\frac{dx}{\sqrt{x(1-x)}}=\arcsin (2x-1),\\
y&=x^{-1}(1-x)^{-(N-1)/4}M(x)^{-1/2}\eta
\end{align}
\end{subequations}
turns the equation
\begin{equation}
\mathcal{L}y=\lambda y+f
\end{equation}
to the standard form
\begin{equation}
-\frac{d^2\eta}{dx^2}+q\eta=\lambda \eta+\hat{f},
\end{equation}
where
\begin{equation}
\hat{f}(x)=x^{-1/2}(1-x)^{-(N-3)/4}M(x)^{-1/2}f(x),
\end{equation}
and
$$q=L_0+\frac{1}{4}\frac{a}{b}
\Big(D\Big(\frac{Da}{a}+\frac{Db}{b}\Big)-\frac{1}{4}
\Big(\frac{Da}{a}+\frac{Db}{b}\Big)^2
+\frac{Da}{a}\Big(\frac{Da}{a}+\frac{Db}{b}\Big)\Big).
$$
Note that $x=0,1/2,1$ are mapped to $\xi=-\pi/2, 0, \pi/2$ and
$$x\sim
\begin{cases}
\displaystyle\frac{1}{4}\frac{1}{(\xi+\pi/2)^2} \quad\mbox{as}\quad x\rightarrow 0 \\
\\
\displaystyle\frac{1}{4}\frac{1}{(\pi/2-\xi)^2}\quad\mbox{as}\quad x\rightarrow 1.
\end{cases}
$$
We see
\begin{equation}
q\sim
\begin{cases}
\displaystyle\frac{2}{(\xi+\pi/2)^2}\quad\mbox{as}\quad \xi\rightarrow -\pi/2 \\
\\
\displaystyle\frac{(N-1)(N-3)}{4}\frac{1}{(\pi/2-\xi)^2}
\quad\mbox{as}\quad \xi\rightarrow \pi/2.
\end{cases}
\end{equation}
Note that $(N-1)(N-3)/4 >3/4$ under the assumption {\bf(B0)}: $N>4$.
Therefore we have

\begin{Proposition}
The operator $\mathfrak{S}_0,\mathsf{D}(\mathfrak{S}_0)=
C_0^{\infty}(-\pi/2,\pi/2), $

\noindent $\mathfrak{S}_0\eta=-d^2\eta/dx^2+q\eta$ in
$L^2(-\pi/2,\pi/2)$ has the Friedrichs extension
$\mathfrak{S}$, a self-adjoint operator,
whose spectrum consists of simple eigenvalues
$\lambda_1<\cdots<\lambda_n<\lambda_{n+1}<\cdots\rightarrow
+\infty$. Thus the operator
$\mathfrak{T}_0, \mathsf{D}(\mathfrak{T}_0)=
C_0^{\infty}(0,1), \mathfrak{T}_0y=\mathcal{L}y$
in

\noindent $L^2(([0,1], x^{3/2}(1-x)^{N/2-1}dx)$ has the self-adjoint
extension $\mathfrak{T}$ with eigenvalues $(\lambda_n)_n$.
\end{Proposition}

Let us fix a positive eigenvalue $\lambda=\lambda_n$ of
$\mathfrak{T}$ and an associated eigenfunction $\Phi$, and
put
\begin{align*}
Y_1(t,x)&=\sin(\sqrt{\lambda}t+\Theta_0)\Phi(x), \\
V_1(t,x)&=\frac{\sqrt{\lambda}}{J(x,0,0)}
\cos(\sqrt{\lambda}t+\Theta_0)\Phi(x),
\end{align*}
$\Theta_0$ being a constant. Then $(Y_1, V_1)$ is a time
periodic solution of the linearized equation
\begin{align*}
&\frac{\partial Y_1}{\partial t}-J(x,0,0)V_1=0, \\
&\frac{\partial V_1}{\partial t}+H_1(x,0,0,0)\mathcal{L}Y_1=0
\end{align*}
( Recall {\bf (B2)}. ) We can claim

\begin{Proposition}
We have
$$
\Phi(x)=
\begin{cases}
C_0+[x]_1\quad\mbox{as}\quad x\rightarrow 0 \\
C_1+[1-x, (1-x)^{N/2}]_1\quad\mbox{as}\quad x\rightarrow 1,
\end{cases}
$$
where $C_0, C_1$ are non-zero constants.
\end{Proposition}

Let us assume\\

{\bf (B4):}  {\it $N/2$ is an integer, that is, $N$ is an even natural number.} \\

Then $\Phi$ is analytic on a neighborhood of $[0,1]$,
and $y^*=\varepsilon Y_1, v^*=\varepsilon V_1$ belongs to
$C^{\infty}(\mathbb{R}\times [0,1])$, $\varepsilon$ being a
parameter. If $|\varepsilon|$ is small, then
$y^*, z^*, v^*, w^* \in U$ for $\forall t \in\mathbb{R}, 0\leq\forall t
\leq 1$. Thus we can apply Theorem 1. The conclusion is

\begin{Corollary}
Assume {\bf (B4)}. There exists $\epsilon_0(T)>0$ such that,
if $|\varepsilon|\leq\epsilon_0(T)$, then we have a solution $(y,v)\in
C^{\infty}([0,T]\times[0,1])$ such that
\begin{equation}
y=\varepsilon Y_1+O(\varepsilon^2),
\quad
v=\varepsilon V_1+O(\varepsilon^2).
\end{equation}
\end{Corollary} 

\textbullet \  Now forget {\bf (B4)} and consider the Cauchy problem
\begin{subequations}
\begin{align}
&\partial_ty-Jv=0,\quad
\partial_tv+H_1\mathcal{L}y+H_2=0, \label{509a} \\
&y|_{t=0}=\psi_0(x),\quad
v|_{t=0}=\psi_1(x), \label{509b}
\end{align}
\end{subequations}
where $\psi_0, \psi_1$ are given functions in
$C^{\infty}([0,1])$. Put
\begin{equation}
y^*=\psi_0(x)+tJ(x,0,0)\psi_1(x),\quad
v^*=\psi_1(x).
\label{510}
\end{equation}
Then we can apply Theorem 1 and get
\begin{Corollary}
There exists a positive number $\delta_0(T)$ such that,
if
$$\max_{k\leq\mathfrak{K}}
\|(D^k\psi_0, D^k\psi_1)\|_{L^{\infty}}\leq \delta_0(T),
$$
then there exists a solution $(y, v)\in
C^{\infty}([0,T]\times[0,1])$ of the Cauchy
problem (\ref{509a})(\ref{509b}).
\end{Corollary}

\section{Open problem}

We hope to consider the following generalized problem.\\

The equations to be considered are of the same type
\begin{subequations}
\begin{align}
&\frac{\partial y}{\partial t}-J(x, y, xDy)v=0, \label{G01a}\\
&\frac{\partial v}{\partial t}+H_1(x, y, xDy, v)\mathcal{L}y+
H_2(x,y, xDy, v, xDv)=0 \label{G01b}
\end{align}
\end{subequations}
with
\begin{equation}
\mathcal{L}=-x(1-x)\frac{d^2}{dx^2}-
\Big(\frac{5}{2}(1-x)-\frac{N}{2}x\Big)\frac{d}{dx}+
L_1(x)\frac{d}{dx}+L_0(x). \label{G02}
\end{equation}

Assume {\bf (B0):} `$N>4$'.

Let us denote by $\mathfrak{A}_{(N)}$ the set of
all smooth functions $f(x)$ of $x \in [0,1[$ such that
$$
f(x)=
\begin{cases}
[x]_0 \quad\mbox{as}\quad x\rightarrow 0 \\
[1-x,(1-x)^{N/2}]_0 \quad\mbox{as}\quad x\rightarrow 1
\end{cases}
$$
and by $\mathfrak{A}_{(N)}^Q(U^p)$,
$U$ being a neighborhood of $0$, the set of
all smooth functions 
$f(x,y_1,\cdots, y_p)$ of
$x \in [0,1[, y_1,\cdots. y_p \in U$ such that 
there are convergent power series
$$\Phi_0(X,Y_1,\cdots, Y_p)=\sum_{k_1+\cdots+ k_p\geq Q}
a^{[0]}_{jk_1\cdots k_p}X^jY_1^{k_1}\cdots Y_p^{k_p}$$
and $$\Phi_1(Z_1, Z_2, Y_1,\cdots, Y_p)=
\sum_{k_1+\cdots+k_p\geq Q}
a^{[1]}_{j_1j_2k_1\cdots k_p}
Z_1^{j_1}Z_2^{j_2}Y_1^{k_1}\cdots Y_p^{k_p}$$
such that
$$f(x,y_1,\cdots, y_p)=\Phi_0(x,y_1,\cdots, y_p) \quad\mbox{for}\quad 0<x\ll 1$$
and 
$$
f(x, y_1,\cdots, y_p)=\Phi_1(1-x, (1-x)^{N/2},y_1,\cdots, y_p)
\quad\mbox{for}\quad 0<1-x\ll 1.$$
Using these notations we assume \\

{\bf (B1'; $N$):} {\it $L_0, L_1\in \mathfrak{A}_{(N)} $ and
$$
L_1(x)=
\begin{cases}
[x]_1\quad\mbox{as}\quad x\rightarrow 0 \\
[1-x,(1-x)^{N/2}]_1\quad\mbox{as}\quad x\rightarrow 1
\end{cases}
$$
and
there is a neighborhood $U$ of $0$ such that
$J\in \mathfrak{A}_{(N)}^0(U^2)$,
$H_1\in\mathfrak{A}_{(N)}^0(U^3)$,
$H_2\in\mathfrak{A}_{(N)}^2(U^4)$.}\\

\textbullet \hspace{5mm} Note that, if $f(x)$ is a function of $\mathfrak{A}_{(N)}$, there is a smooth function
$\Phi \in C^{\infty}([0,1]\times[0,1])$ such that
$$f(x)=\Phi(1-x, (1-x)^{N/2})\qquad \mbox{for}\quad 0\leq\forall x <1, $$
althogh such an analytic two variable function $\Phi$
does not exist.\\

We suppose the assumptions {\bf (B2)}, {\bf (B3)}. 

But in the assumption {\bf (B3)}, ``$f\equiv 0$ as $x\rightarrow 1$" is defined as that there is 
a convergent power series $\Phi(Z_1,Z_2, Y_1,\cdots, Y_p)$ such that
$$f(x,y_1,\cdots, y_p)=(1-x)
\Phi(1-x, (1-x)^{N/2}, y_1, \cdots, y_p)
\quad\mbox{for}\quad 0<1-x\ll 1.
$$
This condition is much stronger than that
$f(1,y_1,\cdots, y_p)=0$ for $\forall y_1\cdots\forall y_p$ if
the assumption {\bf (B4)} :`$N/2$ is an integer' does not hold.
In this sense, we should denote the assumption {\bf (B3)} by
{\bf (B3'; $N$) },writing `$\equiv_{(N)}$' instead of `$\equiv$', to avoid confusions.

Of course if {\bf (B4) }holds, the problem has already been solved. Hence
we are interested in the case that $N/2$ is not an integer,
that is, we assume \\

$\neg${\bf (B4)} : $N$ is not an even integer. \\

So we consider the equations (\ref{G01a})(\ref{G01b}) under the assumptions
{\bf (B0), (B1'; $N$), (B2), (B3'; $N$)} and {\bf $\neg$(B4)}.

Note that the result of the spectral analysis of the linear operator $\mathcal{L}$ is same. In fact we can easily verify that
$$M(x):=\exp\Big[\int_0^x\frac{L_1(x')}{x'(1-x')}dx'\Big]$$
is smooth on $[0,1[$ and enjoys
$$M(x)=
\begin{cases} M_0+[x]_1 \quad\mbox{as}\quad x\rightarrow 0 \\
M_1+[1-x,(1-x)^{N/2}]_1\quad
\mbox{as}\quad x\rightarrow 1
\end{cases}
$$
with positive constants $M_0, M_1$, and therefore, that
eigenfunctions enjoys the claim of Proposition 8.

Hence the problem may be settled as follows:

Fix $T>0$, and functions $y^*, v^* \in \mathfrak{B}_{(N)}([0,T]; U)$.

Here $\mathfrak{B}_{(N)}([0,T]; U)$ stands for the set of all smooth functions $u(t,x)$ of $0\leq t\leq T, 0\leq x<1$, valued in $U$, such that
there is a smooth function 
$\Phi\in C^{\infty}([0,T]\times[0,1]\times[0,1])$ such that
$$u(t,x)=\Phi(t, 1-x, (1-x)^{N/2})$$
for $0\leq t\leq T, 0\leq x <1$. 

Then we seek a solution $(y,v)$ of (\ref{G01a})(\ref{G01b}) of the form
\begin{equation}
y=y^*+\tilde{y}, \qquad v=v^*+\tilde{v} \label{G03}
\end{equation}
such that
\begin{equation}
\tilde{y}|_{t=0}=0, \qquad \tilde{v}|_{t=0}=0. \label{G04}
\end{equation}

In this case
the problem is still open, but is very important in view of
astrophysical applications.

\end{document}